\documentclass{article}

\usepackage{amsmath}
\usepackage{amssymb}
\usepackage{amsthm}
\usepackage{graphicx}
\usepackage{bbm}

\newtheorem{theorem}{Theorem}
\newtheorem{proposition}[theorem]{Proposition}

\begin{document}

\title{The Riff-Shuffle Distribution is Unimodal}

\author{
Stefan Gerhold\thanks{This work was financially supported by the Christian Doppler Research Association (CDG). The author gratefully acknowledges the fruitful collaboration and support by the Bank Austria Creditanstalt (BA-CA) and the Austrian Federal Financing Agency (ÖBFA) through CDG.}\\
Christian Doppler Laboratory for Portfolio Risk Management\\
Vienna University of Technology, Austria\\
{\tt sgerhold at fam.tuwien.ac.at}
}

\maketitle

\begin{abstract}
We show that the probability mass function of the riff-shuffle distribution, also known
as the minimum negative binomial distribution, is unimodal, but in general not log-concave.
\end{abstract}



\vskip 1cm

The riff-shuffle distribution is defined for $0<p<1$ and integers $m \geq 1$ by the probability mass function
\[
  f_k = \binom{m+k-1}{k}(p^m q^k + q^m p^k), \qquad k=0,\dots,m-1,
\]
where $p+q=1$. If cards are taken with probabilities $p$ and $q$ from two decks of $m$ cards, then the number of cards
that have been taken from the remaining deck when one deck has been depleted follows this distribution~\cite{UpBl80}.
Equivalently, it is the distribution of the random variable $X-m$, where $X$ denotes the number of Bernoulli trials with individual success probability $p$ until either
$m$ successes or $m$ failures have been observed.
This explains why it is also called minimum negative binomial distribution~\cite{ZhBuZe00}.
Johnson, Kotz, and Kemp~\cite{JoKeKo05} report that Uppuluri and Blot's~\cite{UpBl80} plots exhibit unimodality of the distribution. In this short note we confirm this observation.

Recall that a sequence $(g_k)_{k=0}^n$ is called log-concave if $g_k^2\geq g_{k-1}g_{k+1}$ holds for
$1\leq k<n$. It is well known that a positive log-concave sequence is unimodal. The riff-shuffle distribution, however,
provides examples of unimodal sequences that are not log-concave.

\begin{proposition}
  The probability mass function of the riff-shuffle distribution is log-concave for $p=\tfrac12$. For
  $p\neq\tfrac12$ and $m$ sufficiently large, it is not log-concave.
\end{proposition}
\begin{proof}
  After cancelling the binomial coefficients, using the relation
  \[
    \binom{u}{v} = \frac{u}{v}\binom{u-1}{v-1},
  \]
  the log-concavity condition $f_k^2\geq f_{k-1}f_{k+1}$ becomes
  \[
    (p^m q^k + q^m p^k)^2 \geq \frac{k(m+k)}{(k+1)(m+k-1)} (p^m q^{k-1} + q^m p^{k-1})(p^m q^{k+1} + q^m p^{k+1}).
  \]
  For $p=\tfrac12$, this reduces to
  \[
    1 \geq \frac{k(m+k)}{(k+1)(m+k-1)},
  \]
  which is easily seen to be true for $k=1,\dots,m-1$. To show the second assertion, we set $k=m-1$. After some simplifications, we see that
  $f_{m-1}^2\geq f_{m-2}f_m$ is equivalent to
  \begin{multline*}
    2(2m-1)(m-1)(pq)^{2m-1}+(m-1)(p^{2m}q^{2m-2}+2(pq)^{2m-1}+p^{2m-2}q^{2m}) \\
    \geq (2m-1)(m-1)\left(\frac{p}{q}+\frac{q}{p}\right)(pq)^{2m-1}.
  \end{multline*}
  The left-hand side of this inequality is $\sim 4m^2(pq)^{2m-1}$ as $m\to\infty$, whereas the right-hand side is $\sim 2(p/q+q/p)m^2(pq)^{2m-1}$.
  Since $p/q+q/p$ is a convex function of $0<p<1$ with minimum $2$ at $p=\tfrac12$, the right-hand side is larger
  than the left-hand side for $p\neq\tfrac12$ and sufficiently large $m$.
\end{proof}

Since the sufficient criterion of log-concavity fails us for $p\neq\frac12$, we now prove directly that the
distribution is unimodal.

\begin{proposition}
  The probability mass function of the riff-shuffle distribution is unimodal, with the maximum
  $\max f_k$ occuring for at most two (adjacent) values of $k$.
\end{proposition}
\begin{proof}
  For $p=\frac12$ it is straightforward to show that $f_0<f_1<\dots<f_{m-2}=f_{m-1}$.
  Now we assume w.l.o.g.\ that $0<p<\tfrac12$.
  Defining the auxiliary function $h(x):=p^m q^x + q^m p^x$, we have
  \[
    \frac{f_{k+1}}{f_k} = \frac{k+m}{k+1} \cdot \frac{h(k+1)}{h(k)}, \qquad k=0,\dots,m-2.
  \]
  Therefore, $f_{k+1}\leq f_k$ is equivalent to $g(k)\geq0$, where
  \[
    g(x):= \frac{h(x)}{h(x+1)} - \frac{x+m}{x+1}.
  \]
  We show that $g$ is concave on $[0,m-2]$. The second derivative
  \[
    \frac{\mathrm{d^2}}{\mathrm{d}x^2} \frac{h(x)}{h(x+1)} =
      \frac{(q-p)(\log p - \log q)^2 p^{m+x} q^{m+x} (p^m q^{x+1}-p^{x+1}q^m)}{(p^m q^{x+1}+p^{x+1}q^m)^3}
  \]
  is non-positive for $0\leq x\leq m-2$, since $(p/q)^{m-x-1}\leq1$ implies
  $p^m q^{x+1}-p^{x+1}q^m \leq 0$. Thus, $h(x)/h(x+1)$ is concave.
  Since $(x+m)/(x+1)$ is a convex function of $x\in[0,m-2]$, the function $g$ is indeed concave.
  At $x=m-2$, the function value
  \begin{align*}
    g(m-2) &= \frac{h(m-2)}{h(m-1)} - \frac{2m-2}{m-1} = \frac{p^m q^{m-2} + q^m p^{m-2}}{p^m q^{m-1} + q^m p^{m-1}} - 2 \\
    &= \frac{1/p^2+1/q^2}{1/p+1/q}-2 > 0
  \end{align*}
  is positive.
  {}From this and concavity we conclude that, if $g(l)\geq0$ for some $l$, then the values $g(k)$, $l<k\leq m-2$, are positive.
  In other words, as soon as there is a descent $f_{l+1}\leq f_l$, the sequence $(f_k)_{k=l+1}^{m-1}$ strictly decreases.
\end{proof}

\bibliographystyle{siam}
\bibliography{../gerhold}

\begin{thebibliography}{1}

\bibitem{JoKeKo05}
{\sc N.~L. Johnson, A.~W. Kemp, and S.~Kotz}, {\em Univariate discrete
  distributions}, Wiley Series in Probability and Statistics,
  Wiley-Interscience [John Wiley \& Sons], Hoboken, NJ, third~ed., 2005.

\bibitem{UpBl80}
{\sc V.~R.~R. Uppuluri and W.~J. Blot}, {\em A probability distribution arising
  in a riff-shuffle}, in Random counts in scientific work, vol 1: {R}andom
  counts in models and structures, Pennsylvania State University Press,
  University Park, 1970, pp.~23--–46.

\bibitem{ZhBuZe00}
{\sc Z.~Zhang, B.~A. Burtness, and D.~Zelterman}, {\em The maximum negative
  binomial distribution}, J. Statist. Plann. Inference, 87 (2000), pp.~1--19.

\end{thebibliography}

\end{document}